\documentclass[reqno]{amsart}
\usepackage{hyperref}

\usepackage[centertags]{amsmath}
\usepackage{amsfonts}
\usepackage{amssymb}
\usepackage{amsthm}
\usepackage{newlfont}
\usepackage{fancyhdr}
\newtheorem{rem}{Remark}

\makeatletter
\@namedef{subjclassname@2010}{%
  \textup{2020} Mathematics Subject Classification}
\makeatother

\fancypagestyle{plain}{\fancyfoot[C]{\vspace*{2\baselineskip}\thepage}}
\theoremstyle{plain}
\newtheorem{theo}{Theorem}[section]
\newtheorem{thmx}{Theorem}[section]

\newtheorem*{conj*}{Denjoy's Conjecture}

\makeatletter
\@namedef{subjclassname@2020}{%
  \textup{2020} Mathematics Subject Classification}
\makeatother

\begin{document}
\title[
\hfil entire solutions of binomial differential equations]
{On a question of Gundersen-Yang concerning entire solutions of binomial differential equations}

\author[ J. R. Long, M. T. Xia, X. X. Xiang \hfil 
]
{Jianren Long, Mengting Xia, Xuxu Xiang}

\address{Jianren Long \newline
School of Mathematical Science, Guizhou  Normal University, Guiyang, 550025,  China.}
\email{longjianren2004@163.com }

\address{Mengting Xia \newline
School of Mathematical Science, Guizhou  Normal University, Guiyang, 550025, China.}
\email{2190331341@qq.com }

\address{Xuxu Xiang \newline
School of Mathematical Science, Guizhou  Normal University, Guiyang, 550025, China.}
\email{1245410002@qq.com }



\subjclass[2020]{34M10, 34M05, 30D35}
\keywords{Nevanlinna theory; Entire solutions; Binomial differential equations; Nonlinear differential equations}

\begin{abstract}
We study the question posed by G. Gundersen and C. C. Yang, in which the following two types of binomial differential equations are investigated,
$$
a(z)f'f''-b(z)(f)^{2}=c(z)e^{2d(z)},~~a(z)ff'-b(z)(f'')^{2}=c(z)e^{2d(z)},
$$
where $a(z)$, $b(z)$ and $c(z)$ are polynomials such that $a(z)b(z)c(z)\not\equiv 0$, $d(z)$ is non-constant polynomial. The explicit forms of entire solutions of the above binomial differential equations are obtained by using the Nevanlinna theory, which gives partial solutions to the question of G. Gundersen and C. C. Yang. In addition, some examples are given to illustrate these results.
\end{abstract}

\maketitle
\numberwithin{equation}{section}
\newtheorem{theorem}{Theorem}[section]
\newtheorem{lemma}[theorem]{Lemma}
\newtheorem{definition}[theorem]{Definition}
\newtheorem{example}[theorem]{Example}
\newtheorem{remark}[theorem]{Remark}
\allowdisplaybreaks

\section{Introduction and main results}

How to characterize the form of meromorphic function $f$ by using the properties of $ff^{(k)}$ is a very interesting topic, such as its zeros, poles and periodicity, and so on. In 1959, Hayman obtained the exact form of $f$ when $f$ and $f''$ only have finitely many zeros and poles, see \cite{hayman} for the details. Langley \cite{langley93} obtained the similar result when $f$ and $f''$ only have finitely many zeros, in which the condition of the number of poles is deleted. In fact, Hayman observed that if $f$ has no zeros, then there exists an entire function $g$ such that $f=\frac{1}{g}$, see \cite[p. 34]{hayman} and \cite[p. 77]{c} for the more details. Thus $f''=\frac{2(g')^{2}-gg''}{g^{3}}$, which leads to the question of when can differential polynomial $2(g')^{2}-gg''$ in an entire function $g$ have no zeros? Many researchers payed attention to this question, at the same time many results have been to obtained, such as \cite{b,fang,k,d,mues71,mues78}, and so on. Recently, Gundersen and Yang \cite{gy} studied these types differential polynomials and the equations associated with them, i.e. they considered the following nonlinear binomial differential equation
\begin{equation}\label{eq1.1}
    ff''-a(z)(f')^2=b(z)e^{2c(z)},
\end{equation}
where $a(z),b(z),c(z)$ are polynomials such that $b(z)\not\equiv 0$ and $c(z)$ is non-constant. For the linear differential equations, a lot of results have been obtained, such as \cite{hitw,g,m,lqt} and therein references.

In what follows, we assume the reader is familiar with the basic notion of Nevanlinna theory, see \cite{c,f} for the details. Let $\rho(f)$ and $\lambda(f)$ denote the order of $f$ and the exponent of convergence of zeros of $f$, respectively.  If $\rho(f)<\infty$, then $S(r,f)=O(\log r)$, if $\rho(f)=\infty$, then $S(r,f)=O(\log rT(r,f))$ for all $r$ outside of an exceptional set with finite logarithmic measure.

We start describe the result from Gundersen and Yang \cite{gy} as follows.

\begin{thmx}\cite{gy}\label{thm1} 
Let f be an entire solution of (\ref{eq1.1})  Then we have one of the following situations:

(1) $b(z)$ is non-constant and $f(z)=p(z)e^{c(z)}$, where p(z) is a polynomial.

(2) $a(z) \equiv \tau$ for some constant $\tau\ne1$, $b(z)\equiv\mu$ for some constant $\mu\ne 0$, $c(z)=\lambda z$ for some constant $\lambda\ne0$, and
$$f(z)=\alpha e^{\lambda z},$$
where $\alpha$ is a constant satisfying $\alpha^2\lambda^2(1-\tau)=\mu$.

(3) $a(z)\equiv1,b(z)\equiv\mu$ for some constant $\mu\ne0,c(z)=\lambda z^2+\gamma z$, and
$$f(z)=\alpha e^{\lambda z^2+\gamma z},$$
where $\alpha,\lambda,\gamma$ are constants such that $2\lambda\alpha^2=\mu$.

(4) $a(z)\equiv1,b(z)\equiv\mu$ for some constant $\mu\ne0,c(z)=\lambda z$ for some constant $\lambda\ne0$, and one of the following three cases occur:
$$f(z)=(\alpha z+\beta) e^{\lambda z},$$
where $\alpha$ and $\beta$ are constants such that $\alpha^2=-\mu$;
$$f(z)=\alpha e^{2\lambda z}+\beta,$$
where $\alpha$ and $\beta$ are constants satisfying $4\alpha\beta\lambda^2=\mu$;
$$f(z)=\alpha e^{\gamma z}+\beta e^{(2\lambda-\gamma)z},$$
where $\alpha,\beta,\gamma$ are constants satisfying
$$\gamma(2\lambda-\gamma)\ne 0 ~ and ~   4\alpha\beta(\lambda-\gamma)^2=\mu.$$
\end{thmx}

In the same paper, Gundersen and Yang posed a analogous problem about the following binomial differential equations:

\begin{equation}\label{eq.1.2}
    a(z)ff''-b(z)(f')^2=c(z) e^{2d(z)},
\end{equation}
\begin{equation}\label{eq.1.3}
    a(z)ff'-b(z)(f'')^2=c(z) e^{2d(z)},
\end{equation}
\begin{equation}\label{eq.1.4}
    a(z)f'f''-b(z)f^2=c(z) e^{2d(z)},
\end{equation}
where $a(z)$, $b(z)$, $c(z)$, $d(z)$ are polynomials such that $a(z)b(z)c(z)\not \equiv 0$ and $d(z)$ is non-constant.

\noindent\textbf{{Question:}}\cite[Question 5]{gy} \emph{What can be said about the entire solutions of the equations \eqref{eq.1.2}-\eqref{eq.1.4}? }

According to the \textbf{Question}, very recently, Yang et. al. \cite{yll} studied the entire solutions of the equations \eqref{eq.1.2} and \eqref{eq.1.4}, the following two results are obtained, respectively.

\begin{thmx}\cite[Theorem 1.1]{yll}\label{thm2}
Assume $\frac{b}{a}$ has at least one simple pole when $\deg (a)\geqslant 1$ and $a+b\not\equiv 0$. If equation \eqref{eq.1.2} admits an entire solution $f$, then $f$ has one of the following expressions:

\begin{itemize}

\item[\textnormal{(i)}]  If $f$ has finitely many zeros, then $$f(z)=p(z)e^{d(z)},$$
where $p$ is a polynomial. In particular,

\begin{itemize}
    \item[\textnormal{(i$_1$)}] if $a+b=\tau$, $c=\mu$, $d=\lambda z$ for some non-zero constants $\tau$, $\mu$, $\lambda$, then $$f(z)=\alpha e^{\lambda z},$$
    where $\alpha$ is a constant satisfying $\tau \lambda^2 \alpha^2=\mu$,

    \item[\textnormal{(i$_2$)}] if $a+b\equiv0$, $a=\tau$, $c=\mu$, $d=\lambda _2z^2+\lambda_1z$ for some non-zero constants $\tau$, $\mu$, $\lambda_1$, $\lambda_2$, then $$f(z)=\alpha e^{\lambda _2z^2+\lambda_1z},$$
    where $\alpha$ is a constant satisfying $2\tau \lambda_2\alpha^2=\mu$,

   \item[\textnormal{(i$_3$)}] if $a+b\equiv0$, $c=\mu a$, $d=\lambda z$ for some non-zero constant  $\mu$, $\lambda$,  then $$f(z)=(\alpha z+\beta)e^{\lambda z},$$
    where $\alpha$, $\beta$ are  constants satisfying $\alpha^2=-\mu$.
\end{itemize}

\item[\textnormal{(ii)}]  If $f$ has infinitely many zeros, but $f'$ has finitely many zeros, then
$$f(z)=T(z)e^{2d}+\beta,$$
where $T\not\equiv0$ is a polynomial satisfying
$$4(a+b)(d')^2T^2+2ad''T^2+4(a+b)d'TT'+aTT''+b(T')^2\equiv0,$$
and $\beta\ne 0$ is a constant. In this case, $\frac{c}{a}$ must be a polynomial. Specifically,

\begin{itemize}
    \item[\textnormal{(ii$_1$)}] if $a+b\equiv0$, then $c=\mu a$, $d=\lambda z$ for some non-zero constants $\mu$, $\lambda$ and $$f(z)=\alpha e^{2\lambda z}+\beta,$$
    where $\alpha$, $\beta$ are  constants satisfying $4\alpha\beta\lambda^2=\mu$,

    \item[\textnormal{(ii$_2$)}] if $a+b\not\equiv0$, then we can deduce $\deg (a)=\deg (b)\geqslant2$, and $\deg (a+b) <\deg (a)$. To be specific, when $k=1$, we have $\deg (a+b)= \deg (a)-2$, $\deg (T)\geqslant1$; and when $k\geqslant2$, we have $\deg (a+b)=\deg (a)-k$.
\end{itemize}

\item[\textnormal{(iii)}]  If both $f$ and $f'$ have infinitely many zeros, then $$f(z)=\alpha e^{\gamma z}+\beta e^{(2\lambda-\gamma)z},$$
where $\alpha$, $\beta$, $\gamma$ are constants satisfying $\gamma(2\lambda-\gamma)\ne 0$ and $4\alpha\beta(\lambda-\gamma)^2=\mu$. In this situation, $a+b\equiv0$, $c=\mu a$, $d=\lambda z$ for some non-zero constants $\mu$, $\lambda$.
\end{itemize}
\end{thmx}

\begin{thmx}\cite[Theorem 1.10]{yll}\label{thm3}
If equation \eqref{eq.1.4} admits an entire solution $f$ satisfying $\lambda(f)<\rho(f)$, where $a$, $b$, $c$, $d$ are polynomials such that $abc\not\equiv0$ and $d$ is non-constant, then $f$ can only be of the form
    $$f(z)=p(z)e^{d(z)},$$
where $p$ is a polynomial.
\end{thmx}


The main purpose of this article is to solve the remaining problem raised by Gundersen and Yang \cite{gy}. We aim to find the exact form of the entire  function solutions of equations \eqref{eq.1.3} and \eqref{eq.1.4}, respectively. The first result show the solution of the \textbf{{Question}} concerning equation \eqref{eq.1.3}.

\begin{theo}\label{the1}
    Let $f$ be an entire solution of equation \eqref{eq.1.3}, where $a$, $b$, $c$, $d$ are polynomials, $abc\not\equiv 0$ and $d$ is non-constant. We have one of the following situations:
\begin{itemize}
    \item [\textnormal{(i)}] If $(ac'-a'c+2acd')f'-acf''\equiv 0$, then
   $$f(z)=\frac{b}{c}(\frac{ac'-a'c}{a^2}+\frac{2cd'}{a})^2e^{2d-c_2}+\frac{e^{c_2}}{c_1},$$
   where $c_1 (\ne 0)$ and $c_2$ are constants. What's more, if $c$ is constant, then $a$, $b$, $d'$ are constants such that $a=8b$ and $(d')^2=c_1$, and $f=\frac{4bc(d')^2}{a^2}e^{2d-c_2}+\frac{e^{c_2}}{c_1}$.
    \item[\textnormal{(ii)}] If $(ac'-a'c+2acd')f'-acf''\not\equiv0$ and $m(r,\frac{f'}{(ac'-a'c+2acd')f'-acf''})=S(r,f)$, then $f$ has one of the following expressions:
    \begin{itemize}
    \item[\textnormal{(ii$_1$)}]  $c$ is non-constant and $f(z)=p(z)e^{d(z)}$, where $p(z)$ and $d(z)$ are polynomials satisfying
   $$app'+ap^2d'-b(p''+2p'd'+pd''+p(d')^2)^2=c.$$

    \item[\textnormal{(ii$_2$)}] $a$, $b$ and $c$ are non-zero constants, and $d=\lambda z$, where $\lambda\ne 0$ is a constant such that $\lambda=-\frac{ad'}{4b(d')^3+a}$, then $f=\delta_2e^{\lambda z}$, where $\delta_2\ne 0$ is a constant such that $-3b(d')^4\delta^2_2=c$.

   \item[\textnormal{(ii$_3$)}] $a$, $b$ and $c$ are non-zero constants, and $d=\lambda z+c_3$, where $\lambda\ne 0$ and $c_3$ are constants, and one of the following three cases occur:
    \begin{itemize}
        \item[\textnormal{(ii$_{31}$)}] $$f(z)=c_1e^{\lambda z},$$ where $c_1$ is non-zero constant such that $c_1^2\lambda(a-b\lambda^3)=ce^{2c_3}$;

        \item[\textnormal{(ii$_{32}$)}]  $$ f(z)=c_1e^{\lambda z}+c_2e^{-2\lambda z},$$
        where $c_1, c_2$ are non-zero constants such that $-9bc_1^2\lambda^4=ce^{2c_3}$;

       \item[\textnormal{(ii$_{33}$)}]  $$f(z)=c_1e^{t_1z}+c_2e^{t_2 z},$$
       where $t_1, t_2$ are two distinct constants such that $t_1^3=t_2^3=\frac{a}{b}$ and $c_1, c_2$ are non-zero constants such that $b c_1c_2t_1^2(t_1-t_2)(t_1+2t_2)=ce^{2c_3}$.
    \end{itemize}

\end{itemize}
    \item[\textnormal{(iii)}]  If $(ac'-a'c+2acd')f'-acf''\not\equiv0$ and $m(r,\frac{f'}{(ac'-a'c+2acd')f'-acf''})\ne S(r,f)$, then the equation \eqref{eq.1.3} does not exist exponential polynomial solutions, where exponential polynomial has the form $f(z)=H_1(z)e^{Q_1(z)}+H_2(z)e^{Q_2(z)}+\cdots +H_k(z)e^{Q_k(z)}$, where $H_j(z)$ and $Q_j(z)$ $(j=1,2,...,k)$ are polynomials in $z$.
\end{itemize}
\end{theo}

\begin{rem}
It is worth noting that the special case $a(z)\equiv1$ is studied by \cite[Theorem 3]{gs}, and our method is different from that of \cite[Theorem 3]{gs}.
\end{rem}

\begin{rem}\label{rema1}
If $abc\equiv0$, then we have the following  situations.

$1)$  $a\equiv0$ and $bc\not \equiv0$. The equation \eqref{eq.1.3} becomes $-b(f'')^2=ce^{2d}$, i.e. $(f'')^2=\frac{-c}{b}e^{2d}$. Since $f$ is an entire function, and $\frac{-c}{b}$ is a polynomial, then $f''$ can only have finitely many zeros. Therefore $f''=Re^{d}$, where $R$ is a polynomial.

$2)$ $b\equiv0$ and $ac\not\equiv0$. The equation \eqref{eq.1.3} becomes $aff'=ce^{2d}$, i.e. $ff'=\frac{c}{a}e^{2d}$. Since $\frac{c}{a}$ is a polynomial, and $ff'=f^2\frac{f'}{f}=\frac{c}{a}e^{2d}$, then $f$ can only have finitely many zeros, and  $T(r,\frac{f'}{f})=N(r,\frac{f'}{f})+m(r,\frac{f'}{f})=S(r,f)$. Since $T(r,f^2)\leqslant T(r, \frac{f}{f'}\frac{c}{a}e^{2d})$, then $\rho(f)<\infty$. Since  $(\frac{f}{e^d})^2=\frac{f}{f'}\frac{c}{a}$, then $2T(r,\frac{f}{e^d})=S(r,f)$. Therefore, $f=Re^d$, where $R$ is a polynomial.

$3)$ $c\equiv0$ and $ab\not\equiv0$. The equation \eqref{eq.1.3} becomes $aff'=b(f'')^2$, i.e. $a\frac{f'}{f}=b(\frac{f''}{f})^2$. We claim that $f$ can only have finitely many zeros, otherwise if $z_0$ is zero with multiplicity of $t$ of $f$ (not the zero of $a$ and $b$), then we have $t+t-1=2t-4$, which is impossible. From Wiman-Valiron theory \cite[Theorem 3.2]{g}, we have $|a|\frac{V}{|z|}=|b|\frac{V^4}{|z|^4}(1+o(1))$ for $|z|=r\notin E$, where $V$ is the central index of $f$, $E$ is a set with finite logarithmic measure. Then by \cite[Lemma 1.1.2]{g}, for sufficiently large $r$, $V^3=\frac{|a|}{|b|}|z|^3$. Let $a(z)=a_nz^n+\cdots +a_0$, $b(z)=b_mz^m+ \cdots +b_0$, where $a_n$ and $b_m$ are non-zero constants. Then for sufficiently large $r$, $\frac{1}{2}|a_n|r^n\leqslant |a| \leqslant 2|a_n|r^n$, and $\frac{1}{2}|b_m|r^m\leqslant |b| \leqslant 2|b_m|r^m$. Therefore, $\frac{1}{4}|\frac{a_n}{b_m}|r^{n+3-m}\leqslant V^3\leqslant 4|\frac{a_n}{b_m}|r^{n+3-m}$. If $n+3\leqslant m$, then it implies that $V$ is a constant, this is a contradiction. Thus $n+3>m$ and $V=O(r^{n+3-m})$, which implies that $f$ has finite order. Hence $f=Re^P$, where $R$ and $P$ are polynomials.
\end{rem}

The following result shows the solution of the \textbf{{Question}} concerning the equation \eqref{eq.1.4}, in which our condition is relax than Theorem \ref{thm3}.

\begin{theo}\label{the2}
    Let f be an entire solution of equation \eqref{eq.1.4}, where $a$, $b$, $c$, $d$ are polynomials such that  $abc\not\equiv 0$ and $d$ is non-constant. Then we have one of the following situations:
\begin{itemize}
    \item[\textnormal{(i)}] There exists at least one of $a, b, c, d'$ that is not constant, then we have one of the following three cases occur:
\begin{itemize}
    \item[\textnormal{(i$_1$)}]  $c$ is non-constant and $f(z)=P(z)e^{d(z)}$, where $P(z)$ is a polynomial satisfying
   $$(aP'+aPd')[P''+2P'd'+Pd''+P(d')^2]-bP^2=c.$$
   \item[\textnormal{(i$_2$)}] $c$ is a non-zero constant, one of the following two cases occur:
    \begin{itemize}
          \item[\textnormal{(i$_{21}$)}] $d'=\lambda$, where $\lambda$ is a non-zero constant and at least one of $a$ and $b$ is not constant, then $\deg (a)=\deg (b)$, and $f(z)=Pe^{d(z)}$, where $P$ satisfying
          $$aP'P''+2\lambda a(P')^2+3a\lambda ^2PP'+\lambda aPP''+P^2(\lambda^3a-b)\equiv c.$$

         \item[\textnormal{(i$_{22}$)}] $d'$ is non-constant and $b-ad'd''-a(d')^3=M$, where $M$ is a non-zero polynomial, one of the following cases occur:

         ${\textnormal{1)}}$ $\deg (M)=0$, and
          $$f(z)=Pe^{d(z)},$$
          where $P$ is a constant satisfying $P^2M=c$;

         ${\textnormal{2)}}$ $\deg (M)\geqslant 1$, $\deg(M)=\deg(a)+2\deg (d)-3$, and $$f(z)=P(z)e^{d(z)},$$
         where $P(z)$ is a polynomial satisfying $$aP'P''+2a(P')^2d'+aPP'd''+3aPP'(d')^2+aPP''d'+P^2M=c.$$
    \end{itemize}
   \item[\textnormal{(i$_3$)}]  $\lambda(f)\leqslant\rho(f)<\infty.$
\end{itemize}
    \item[\textnormal{(ii)}] All of $a, b, c, d'$ are non-zero constants. Let $d=\lambda z+c_3$, where $\lambda\ne 0$ and $c_3$ are constants, then we have one of the following three cases:
\begin{itemize}
    \item[\textnormal{(ii$_1$)}]  $f(z)=c_1e^{\lambda z}$, where $c_1\ne 0$ satisfies $c_1^2(a\lambda^3-b)=ce^{2c_3}$;

   \item[\textnormal{(ii$_2$)}]  $$ f(z)=c_1e^{\lambda z}+c_2e^{-\frac{1}{2}\lambda z},$$
   where $c_1$, $c_2$ are non-zero constants such that $\frac{9}{8}c_1^2a\lambda^3=ce^{2c_3}$;

   \item[\textnormal{(ii$_3$)}]  $$f(z)=c_1e^{t_1z}+c_2e^{t_2 z},$$
   where $t_1$, $t_2$ are two distinct constants such that $t_1^3=t_2^3=\frac{b}{a}$ and $c_1, c_2$ are non-zero constants such that $c_1c_2t_1(t_2-t_1)(t_2+2t_1)=ce^{2c_3}$.
\end{itemize}
\end{itemize}
\end{theo}

\begin{rem}
 Theorem \ref{the2} improves \cite[Theorem 1.5]{wlh} and Theorem \ref{thm3}. On the one hand, it is worth noting that when $a(z)\equiv 1$, Theorem \ref{the2} reduce to the case in \cite[Theorem 1.5]{wlh}, but our theorem remove the restriction condition:  $f'$ only has simple zero in \cite[Theorem 1.5]{wlh}. On the other hand, the restriction $\lambda(f)<\rho(f)$ in Theorem \ref{thm3} is removed in Theorem \ref{the2}, and our theorem gives more solutions of equation \eqref{eq.1.4}.
\end{rem}

\begin{rem}
If $abc\equiv0$, then we have the following situations.

$1)$ $a\equiv0$ and $bc\not\equiv0$. The equation \eqref{eq.1.4} becomes $-bf^2=ce^{2d}$, i.e. $f^2=\frac{-c}{b}e^{2d}$. Since $f$ is an entire function, then $\frac{-c}{b}$ is a polynomial, which implies that $f$ can only have finitely many zeros. Thus $f=Re^{d}$, where $R^2=\frac{-c}{b}$.

$2)$ $b\equiv0$ and $ac\not\equiv0$. The equation \eqref{eq.1.4} becomes $af'f''=ce^{2d}$, i.e. $f'f''=\frac{c}{a}e^{2d}$. Since $f$ is an entire function, $\frac{c}{a}$ is a polynomial, and  $f'f''=(f')^2\frac{f''}{f'}=\frac{c}{a}e^{2d}$, then $f'$ can only have finitely many zeros. This implies  $T(r,\frac{f''}{f'})=N(r,\frac{f''}{f'})+m(r,\frac{f''}{f'})=S(r,f)$. Since $T(r,(f')^2)\leqslant T(r, \frac{f'}{f''}\frac{c}{a}e^{2d})$, then $\rho(f')<\infty$. Since $(\frac{f'}{e^d})^2=\frac{f'}{f''}\frac{c}{a}$, then $2T(r,\frac{f'}{e^d})=S(r,f)$. Therefore, $f'=Re^d$, where $R$ is a polynomial.

$3)$ $c\equiv0$ and $ab\not\equiv0$. The equation \eqref{eq.1.4} becomes $af'f''=bf^2$, i.e. $a\frac{f'}{f}\frac{f''}{f}=b$. We claim that $f$ can only have finitely many zeros, otherwise if $z_0$ is zero with multiplicity of $k$ of $f$ (not the zero of $a$ and $b$), then we have $k-1+k-2=2k$, which is impossible. From Wiman-Valiron theory \cite[Theorem 3.2]{g}, we have $|a|\frac{V^3}{|z|^3}=|b|(1+o(1))$ for $|z|=r\notin E$, where $V$ is the central index of $f$, $E$ is a set with finite logarithmic measure. Then by \cite[Lemma 1.1.2]{g}, for sufficiently large $r$, $V^3=\frac{|b|}{|a|}|z|^3$. Let $a(z)=a_nz^n+\cdots +a_0$, $b(z)=b_mz^m+ \cdots +b_0$, where $a_n$ and $b_m$ are non-zero constants. Then for sufficiently large $r$, $\frac{1}{2}|a_n|r^n\leqslant |a| \leqslant 2|a_n|r^n$, and $\frac{1}{2}|b_m|r^m\leqslant |b| \leqslant 2|b_m|r^m$. Therefore, $\frac{1}{4}|\frac{b_m}{a_n}|r^{m+3-n}\leqslant V^3\leqslant 4|\frac{b_m}{a_n}|r^{m+3-n}$. If $m+3\leqslant n$, then it implies that $V$ is a constant, this is a contradiction. Thus $m+3>n$ and $V=O(r^{m+3-n})$, which implies that $f$ has finite order. Hence $f=Re^P$, where $R$ and $P$ are polynomials.
\end{rem}

\section{Some Examples}

In this section, some examples are given to illustrate our results. Examples \ref{exa6}-\ref{exa8} will to show the partial cases of Theorem \ref{the1}.

\begin{example}\label{exa6}
    Let $a(z)=16, b(z)=2, c(z)=64, d(z)=z$. Then
$$f(z)=e^{2z}+2$$
satisfies the equation $16f(z)f'(z)-2(f''(z))^2=64e^{2z}$.
This is the situation ${\textnormal{(i)}}$ of Theorem \ref{the1}.
\end{example}

\begin{example}\label{exa7}
    Let $a(z)=28, b(z)=1, c(z)=12, d(z)=3z$. Then
$$f(z)=2e^{3z}$$
satisfies the equation $28f(z)f'(z)-(f''(z))^2=12e^{6z}$.
This is the situation ${\textnormal{(ii$_{31}$)}}$ of Theorem \ref{the1}.
\end{example}

\begin{example}\label{exa8}
    Let $a(z)=-128, b(z)=2, c(z)=-288, d(z)=2z$. Then
$$f(z)=e^{2z}+2e^{-4z}$$
satisfies the equation $-128f(z)f'(z)-2(f''(z))^2=-288e^{4z}$.
This is the situation ${\textnormal{(ii$_{32}$)}}$ of Theorem \ref{the1}.
\end{example}

The following Examples \ref{exa1}-\ref{exa5} are used to show the partial cases of Theorem \ref{the2}.

\begin{example}\label{exa1}
    Let $a(z)=3z+4, b(z)=8(3z+4), c(z)=4(3z+4)^2, d(z)=2z$. Then
$$f(z)=(z+1)e^{2z}$$
satisfies the equation $(3z+4)f'(z)f''(z)-8(3z+4)(f(z))^2=4(3z+4)^2e^{4z}$.
This is the situation ${\textnormal{(i$_{1}$)}}$ of Theorem \ref{the2}.
\end{example}

\begin{example}\label{exa2}
     Let $a(z)=1, b(z)=2z, c(z)=8z^3+2z, d(z)=z^2$. Then
$$f(z)=e^{z^2}$$
satisfies the equation $f'(z)f''(z)-2z(f(z))^2=(8z^3+2z)e^{2z^2}$.
This is the situation ${\textnormal{(i$_1$)}}$ of Theorem \ref{the2}.
\end{example}

\begin{example}\label{exa3}
    Let $a(z)=2, b(z)=16z^3+24z^2+20z+7, c(z)=-9, d(z)=z^2+z$. Then
$$f(z)=3e^{z^2+z}$$
satisfies the equation $2f'(z)f''(z)-(16z^3+24z^2+20z+7)(f(z))^2=-9e^{2z^2+2z}$.
This is the situation ${\textnormal{(i$_{22}$)}}-1)$ of Theorem \ref{the2}.
\end{example}

\begin{example}\label{exa4}
    Let $a(z)=3, b(z)=1, c(z)=88, d(z)=2z$. Then
$$f(z)=2e^{2z}$$
satisfies the equation $3f'(z)f''(z)-2(f(z))^2=88e^{4z}$.
This is the situation ${\textnormal{(ii$_{1}$)}}$ of Theorem \ref{the2}.
\end{example}

\begin{example} \cite{gy}\label{exa5}
     Let $\alpha,\beta,\lambda,\mu,$  are non-zero constants such that $\lambda^3=\mu ^3$ and $\lambda \ne \mu$ , $a(z)=1, b(z)=\lambda^3, c(z)=\alpha\beta\lambda(\mu -\lambda )(2\lambda +\mu), d(z)=\frac{(\lambda +\mu )z}{2} $. Then
$$f(z)=\alpha e^{\lambda z}+\beta e^{\mu z}$$
satisfies the equation $f'(z)f''(z)-\lambda ^3(f(z))^2=\alpha\beta\lambda(\mu -\lambda )(2\lambda +\mu )e^{(\lambda +\mu )z}$.
This is the situation ${\textnormal{(ii$_{3}$)}}$ of Theorem \ref{the2}.
\end{example}


\section{Proof of Theorem \ref{the1}}









Let $f$ be an entire solution of \eqref{eq.1.3}. Taking the logarithm of both sides of the equation \eqref{eq.1.3}, and then take the derivative,
$$\frac{a'ff'+a(f')^2+aff''-b'(f'')^2-2bf''f'''}{aff'-b(f'')^2}=\frac{c'+2cd'}{c}.$$
Then we have
\begin{equation}
\label{3.1}
    f'(acf'-hf)=f''(2bcf'''+sf''-acf),
\end{equation}
where $h=ac'-a'c+2acd'$, $s=b'c-bc'-2bcd'$. We claim that $h\not\equiv0$, otherwise, $ac'-a'c+2acd'\equiv0$. Since $ac\not\equiv0$, divide both sides of $h$ by $ac$, we have $\frac{c'}{c}-\frac{a'}{a}=-2d'$, then $e^{-2d}=\delta_1\frac{c}{a}$, $\delta_1\ne0$ is a constant, this is impossible. Using the same proof for $h\not\equiv 0$, we have  $s\not\equiv0$. Rewriting \eqref{3.1} as
\begin{eqnarray*}
      ac(f')^2-2bcf''f'''-s(f'')^2=f(hf'-acf'').
\end{eqnarray*}
Divide both sides of the above equation by $(f')^2$ to get
\begin{equation}
\label{3.2}
    ac-\frac{2bcf''f'''+s(f'')^2}{(f')^2}=\frac{f}{f'}\frac{hf'-acf''}{f'}.
\end{equation}

\indent\textbf{Case 1.} $hf'-acf''\equiv0$. We can see from the equation \eqref{eq.1.3} that $f'$ and $f''$ have only finitely many common zeros. Since $hf'-acf''\equiv0$, then $f'$ and $f''$ have finitely many zeros and $f'=pe^{l}$, where $p$ is a polynomial and $l$ is a non-zero entire function.  Since $f''=(p'+pl')e^l$, then $l$ is a polynomial. Substituting $f'$ and $f''$ into equation \eqref{eq.1.3}, we have
\begin{equation}
\label{3.3}
    f=\frac{b}{ap}(p'+pl')^2e^l+\frac{c}{ap}e^{2d-l}.
\end{equation}
Substituting $f'$ and $f''$ into $hf'-acf''\equiv0$, we have
\begin{equation}
\label{3.4}
    apc'-acp'-a'cp+acp(2d'-l')=0.
\end{equation}
Noting that $acp\not \equiv0$, by \eqref{3.4} we deduce that $\frac{ap}{c}$ and $2d-l$ are both constants. Set $\frac{ap}{c}=c_1$, $2d-l=c_2$, where $c_1 (\ne 0)$ and $c_2$ are constants. then \eqref{3.3} becomes
\begin{equation}
\label{3.5}
    f=\frac{b}{c}(\frac{ac'-a'c}{a^2}+\frac{2cd'}{a})^2e^{2d-c_2}+\frac{e^{c_2}}{c_1}.
\end{equation}
What's more, if $c$ is constant, $ap$ must be constant, that is, $a$ and $p$ must be constants. By \eqref{3.3} and $f'=pe^{l(z)}$, we have
\begin{eqnarray*}
     (\frac{b}{ap}(p'+pl')^2)'+\frac{b}{ap}(p'+pl')^2l'=p.
\end{eqnarray*}
Since $a$ and $p$ are constants, from the above equation we have both $b$ and $l'=2d'$ have to be constant. Then equation \eqref{3.5} becomes $$f=\frac{4bc(d')^2}{a^2}e^{2d-c_2}+\frac{e^{c_2}}{c_1},$$
substituting it into equation \eqref{eq.1.3}, we obtain
$$(\frac{32b^2c^2(d')^4}{a^3}-\frac{256b^3c^2(d')^4}{a^4})e^{4d-2c_2}+(\frac{8bc(d')^2}{ac_1}-c)e^{2d}=0.$$
Applying \cite[Theorem 1.51]{f} to the above equation, we have $\frac{32b^2c^2(d')^4}{a^3}-\frac{256b^3c^2(d')^4}{a^4}=0$, i.e. $a=8b$ and $\frac{8bc(d')^2}{ac_1}-c=0$, i.e. $(d')^2=c_1$. $f$ has the form of $f=\frac{4bc(d')^2}{a^2}e^{2d-c_2}+\frac{e^{c_2}}{c_1}$. This is the conclusion ${\textnormal{(i)}}$ of Theorem \ref{the1}.

\indent\textbf{Case 2.} $hf'-acf''\not\equiv0$ and $m(r,\frac{f'}{hf'-acf''})=S(r,f)$. Applying \cite[Theorem 1.22]{f} to \eqref{3.2}, we have
\begin{equation}
\label{3.6}
    m(r,\frac{f}{f'})=S(r,f).
\end{equation}
From \eqref{3.1}, we set
\begin{equation}
\label{3.7}
    \gamma=\frac{2bcf'''+sf''-acf}{f'}=\frac{acf'-hf}{f''}.
\end{equation}

Looking at equation \eqref{eq.1.3}, we find that $f'$ and $f''$ have only finitely many common zeros. By \eqref{3.1}, we get $N(r,\gamma)=O(\log r)$. By \eqref{3.6} and \eqref{3.7}, we have $m(r,\gamma)=S(r,f)$. Hence, $T(r,\gamma)=S(r,f).$
Equation \eqref{3.7} gives two linear differential equations:
\begin{equation}
\label{3.8}
    \gamma f''-acf'+hf=0,
\end{equation}
\begin{equation}
\label{3.9}
    2bcf'''+sf''-\gamma f'-acf=0.
\end{equation}
By differentiating \eqref{3.8}, we have
\begin{equation}
\label{3.10}
    \gamma f'''+(\gamma' -ac)f''+(h-a'c-ac')f'+h'f=0.
\end{equation}
From \eqref{3.9} and \eqref{3.10}, we get
\begin{equation}
\label{3.11}
    (2bc\gamma'-2abc^2-s\gamma)f''+(2bch-2a'bc^2-2abcc'+\gamma^2)f'+(2bch'+ac\gamma)f=0.
\end{equation}
From \eqref{3.8} and \eqref{3.11}, we get
\begin{equation}
\label{3.12}
    uf'=vf,
\end{equation}
where
\begin{eqnarray}
    u&=&\gamma^3+(2bch-2a'bc^2-2abcc'-acs)\gamma+2abc^2\gamma'-2a^2bc^3, \label{3.13} \\
    v&=&-ac\gamma^2-(2bch'+hs)\gamma+2bch\gamma'-2abc^2h.   \label{3.14}
\end{eqnarray}

Since $ff'\not\equiv0$, then we can divide the discussion into two cases: $uv\not\equiv 0$ or $u\equiv v\equiv0$.

\indent\textbf{Subcase 2.1.}  $uv\not\equiv 0$. We claim that $N(r,\frac{1}{f})=S(r,f)$, otherwise, assume $z_0$ is a zero of $f$ (not the zero of $a$, $b$, $c$, $u$, $\gamma$, $h$), from \eqref{3.12} we get $f'(z_0)=0$, and from \eqref{3.8} we get $f''(z_0)=0$. Substituting $z_0$ into \eqref{eq.1.3}, the left side of \eqref{eq.1.3} is zero and the right side is not zero, this is impossible. Thus $N(r,\frac{1}{f})=S(r,f)$. From \eqref{eq.1.3},
$$\frac{1}{f^2}=\frac{e^{-2d}}{c}(a\frac{f'}{f}-b(\frac{f''}{f})^2).$$
From the first fundamental theorem, the above formula and \cite[Theorem 1.22]{f}, we get
\begin{eqnarray*}
    2T(r,f) = T(r,f^2)&=&m(r,\frac{1}{f^2})+N(r,\frac{1}{f^2})+O(1)  \\
    &\leqslant&m(r,\frac{e^{-2d}}{c}(a\frac{f'}{f}-b(\frac{f''}{f})^2))+S(r,f)  \\
    &\leqslant&m(r,e^{-2d})+m(r,\frac{f'}{f})+2m(r,\frac{f''}{f})+S(r,f)  \\
    &\leqslant&O(r^k)+S(r,f),
\end{eqnarray*}
where $k=\deg (d) >0$. Thus, $\rho(f)\leqslant k$. Moreover, from \eqref{eq.1.3} we have
\begin{eqnarray*}
    T(r,e^{2d})=m(r,e^{2d})&=&m(r,\frac{aff'-b(f'')^2}{c})  \\
    &\leqslant&m(r,\frac{aff'-b(f'')^2}{cf^2})+m(r,f^2)+S(r,f)  \\
    &\leqslant&2T(r,f)+S(r,f).
\end{eqnarray*}
Then we have $\rho (f)\geqslant k$. Hence, $\rho (f)=k$. Since $f$ has finite order,  then $N(r,\frac{1}{f})=O(\log r)$. This shows that $f$ has only finitely many zeros, then $f(z)=p_1(z)e^{w(z)}$,  $p_1(z)$ is a polynomial and $w(z)$ is a polynomial of degree $k$. Substituting $f(z)=p_1(z)e^{w(z)}$ into equation \eqref{eq.1.3}, we obtain
$$[ap_1p_1'+ap_1^2w'-b(p_1''+2p_1'w'+p_1w''+p_1(w')^2)^2]e^{2w}=ce^{2d}.$$
There's only one constant difference between $w(z)$ and $d(z)$. Therefore, we can write $f(z)=p(z)e^{d(z)}$, where $p(z)$ is a polynomial such that
\begin{equation*}
    app'+ap^2d'-b(p''+2p'd'+pd''+p(d')^2)^2= c,
\end{equation*}
where $c$ is a non-constant polynomial. This is the conclusion ${\textnormal{(ii$_1$)}}$ of Theorem \ref{the1}.

\indent\textbf{Subcase 2.2.}  $u\equiv v\equiv 0$.

From $u\equiv 0$ and \eqref{3.13}, we get $\gamma\ne 0$. First we prove $\gamma$ is an entire function, otherwise suppose that $\gamma$ has a pole $z_1$ with multiplicity of $m$, and from \eqref{3.13} knows that $z_1$ is a pole of $u$ with multiplicity of $3m$, but $u\equiv 0$, which is a contradiction. By \eqref{3.13} and \eqref{3.14} we have $hu-acv\equiv0$, i.e.
$$h\gamma^2+a^2c^2\gamma+2abc^2h'+2bch^2-2a'bc^2h-2abcc'h\equiv0.$$

 Next we show $\gamma $ is a polynomial. Otherwise,
 suppose $\gamma$  is transcendental. Since $h\not\equiv0$, using Clunie Lemma \cite[Lemma 3.3]{c} to the above equation, then $T(r, \gamma)=m(r,\gamma)=S(r,\gamma)$, which is impossible.

If there exists at least one of $s$ and $a$ is a non-constant polynomial, then by \cite[Remark 1 of Theorem 4.1]{g}, the entire solution of the equation \eqref{3.11} is of finite order. Thus $\lambda(f)\leqslant \rho(f)<\infty$.

If $s$ and $a$ are both constants $(\ne 0)$, $s=b'c-bc'-2bcd'$, $d'\ne 0$, then the term $2bcd'$ must be constant. Hence, $a$, $b$, $c$, $d'$ and $s$ are non-zero constants, and we know from the equation \eqref{3.13} that $2bch-2a'bc^2-2abcc'-acs$, $2abc^2$ and $-2a^2bc^3$ are constants, then we have $\gamma $ is also constant. We write  \eqref{3.11} as
\begin{equation}
\label{3.15}
    P_2f''+P_1f'+P_0f=0,
\end{equation}
where $P_2=2bc\gamma'-2abc^2-s\gamma$, $P_1=2bch-2a'bc^2-2abcc'+\gamma^2$, $P_0=2bch'+ac\gamma$. Since $a$, $b$, $c$, $d'$, $s$ and $\gamma$ are non-zero constants, then $P_2=2bcd'\gamma-2abc^2$, $P_1=4abc^2d'+\gamma^2$, $P_0=ac\gamma(\ne 0)$.

If $P_2=0$, i.e. $\gamma=\frac{ac}{d'}$, the above equation becomes $P_1f'+P_0f=0$. Since $P_0\ne 0$, we have $P_1\ne 0$, then $f=\delta_2 e^{-\frac{P_0}{P_1}z}$, where $\delta_2\ne 0$ is a constant, i.e. $f=\delta_2e^{-\frac{ad'}{4b(d')^3+a}z}$, substituting it into \eqref{eq.1.3}, we obtain
\begin{eqnarray*}
    -\delta^2_2(\frac{a^2d'}{4b(d')^3+a}+\frac{a^4b(d')^4}{(4b(d')^3+a)^4})e^{-\frac{2ad'}{4b(d')^3+a}z}=ce^{2d}.
\end{eqnarray*}
Then we get $d=-\frac{ad'}{4b(d')^3+a}z$, i.e. $d'=-\frac{ad'}{4b(d')^3+a}$, then we have $a=-2b(d')^3$. From the above equation we get $-3b(d')^4\delta^2_2=c$. This is the conclusion ${\textnormal{(ii$_{2}$)}}$ of Theorem \ref{the1}.

If $P_2\ne 0$, then assume that the characteristic equation of \eqref{3.15} is
$$t^2+\frac{P_1}{P_2}t+\frac{P_0}{P_2}=0,$$
where $t_1$ and $t_2$ are two characteristic roots. Since $P_0\ne0$, then $t_1$ and $t_2$ are non-zero constants.

\indent \textbf{Subcase 2.2.1.} $t_1\ne t_2$.
Then
\begin{equation}
\label{3.16}
    f(z)=c_1e^{t_1z}+c_2e^{t_2z},
\end{equation}
where $c_1$ and $c_2$ are constants. Next we will discuss the case of whether $c_1$, $c_2$ is zero.

\indent \textbf{Subcase 2.2.1.1.} $c_1c_2=0$. Without loss of generality, assume that $c_2=0$.  Since $d'$ is a constant, then assume $d=\lambda z+c_3$, where $\lambda (\ne 0)$ and $c_3$ are constants. Substituting $f(z)=c_1e^{t_1z}$ and  $d=\lambda z+c_3$ into equation \eqref{eq.1.3} to get
$$c_1^2t_1(a-bt_1^3)e^{2t_1z}=ce^{2c_3}e^{2\lambda z},$$
which implies that $\lambda=t_1$, and $c_1^2\lambda(a-b\lambda^3)=ce^{2c_3}$. Then the solution $f$ is of form $f(z)=c_1e^{\lambda z}$. This is the conclusion ${\textnormal{(ii$_{31}$)}}$ of Theorem \ref{the1}.

\indent \textbf{Subcase 2.2.1.2.}  $c_1c_2\ne0$. Substituting \eqref{3.16} and $d=\lambda z+c_3$ into equation \eqref{eq.1.3}, we have
\begin{equation}
\label{3.17}
    c_1^2t_1(a-b t_1^3)e^{2t_1z}+c_2^2t_2(a -b t_2^3)e^{2t_2z}+c_1c_2(a t_1+a t_2-2b t_1^2t_2^2)e^{(t_1+t_2)z}=ce^{2c_3}e^{2\lambda z}.
\end{equation}
Applying \cite[Theorem 1.51]{f} to equation \eqref{3.17}, then $\lambda= t_1$ or $\lambda=t_2$ or $\lambda=\frac{t_1+t_2}{2}$.

If $\lambda=t_1$, then \eqref{3.17} can be written as
$$ [c_1^2t_1(a-b t_1^3)-ce^{2c_3}]e^{2t_1z}+c_2^2t_2(a -b t_2^3)e^{2t_2z}+c_1c_2(a t_1+a t_2-2bt_1^2t_2^2)e^{(t_1+t_2)z}=0.$$
By applying \cite[Theorem 1.51]{f} to the above equation, we have
$$c_1^2t_1(a-b t_1^3)-ce^{2c_3}=c_2^2t_2(a -b t_2^3)=c_1c_2(a t_1+a t_2-2bt_1^2t_2^2)=0,$$
which gives $a=b t_2^3$, and $bc_1c_2t_2^2(t_2-t_1)(t_2+2t_1)=0$. Since $t_1\ne t_2$, we have $\lambda=t_1=-\frac{t_2}{2}$. Combining this with the equation  \eqref{3.16}, we have $f(z)=c_1e^{\lambda z}+c_2e^{-2\lambda z}$ and $-9b c_1^2
\lambda^4=ce^{2c_3}$. This is the conclusion ${\textnormal{(ii$_{32}$)}}$ of Theorem \ref{the1}.

If $\lambda=t_2$, using the same method as $\lambda=t_1$ we get  $\lambda=t_2=-\frac{t_1}{2}$. Then from \eqref{3.16},  we have $f(z)=c_1e^{-2\lambda z}+c_2e^{\lambda z}$ and $-9b c_2^2
\lambda^4=ce^{2c_3}$. This is the conclusion ${\textnormal{(ii$_{32}$)}}$ of Theorem \ref{the1}.

If $\lambda=\frac{t_1+t_2}{2}$, then \eqref{3.17} becomes
$$ c_1^2t_1(a-bt_1^3)e^{2t_1z}+c_2^2t_2(a-bt_2^3)e^{2t_2z}+[c_1c_2(at_1+at_2 -2b t_1^2t_2^2)-ce^{2c_3}]e^{(t_1+t_2)z}=0.$$
By applying \cite[Theorem 1.51]{f} to the above equation, we have
$$c_1^2t_1(a-bt_1^3)=c_2^2t_2(a-bt_2^3)=c_1c_2(at_1+at_2 -2b t_1^2t_2^2)-ce^{2c_3}=0.$$
We get $t_1^3=t_2^3=\frac{a}{b}$ and $bc_1c_2t_1^2(t_1-t_2)(t_1+2t_2)=ce^{2c_3}$. Then from \eqref{3.16}, we have $f(z)=c_1e^{t_1z}+c_2e^{t_2z}$. This is the conclusion ${\textnormal{(ii$_{33}$)}}$ of Theorem \ref{the1}.

\indent\textbf{Subcase 2.2.2.}  $t_1=t_2$. Then
\begin{equation}
\label{3.18}
    f(z)=(c_4+c_5 z)e^{t_1 z},
\end{equation}
where $c_4, c_5$ are constants. Substituting \eqref{3.18} into equation \eqref{eq.1.3}, we get
$$[a(c_4+c_5z)(c_5+c_4t_1+c_5t_1z)-b(2c_5t_1+c_4t_1^2+c_5t_1^2z)^2]e^{2t_1z}=ce^{2c_3} e^{2\lambda z},$$
which implies that $\lambda=t_1$ and
\begin{equation}
\label{3.19}
    a(c_4+c_5z)(c_5+c_4t_1+c_5t_1z)-b(2c_5t_1+c_4t_1^2+c_5t_1^2z)^2=ce^{2c_3}.
\end{equation}
By considering the coefficients at $z^2$ and $z$, respectively, we have
$$c_5^2t_1(a-bt_1^3)=0,$$
and$$ac_5^2+2ac_4c_5t_1-4bc_5^2t_1^3-2bc_4c_5t_1^4=0.$$
If $c_5\ne 0$, then $a=bt_1^3$, and $-3bc_5^2t_1^3=0$, which implies that $c_5=0$, this is a contradiction. Thus $c_5=0$. Substituting $c_5=0$ into \eqref{3.19}, we obtain $c_4^2\lambda(a-b\lambda^3)=ce^{2c_3}$. Then the solution $f$ is of form $f(z)=c_4e^{\lambda z}$. This is the conclusion ${\textnormal{(ii$_{31}$)}}$ of Theorem \ref{the1}.

\indent\textbf{Case 3.} $hf'-acf''\not\equiv0$ and $m(r,\frac{f'}{hf'-acf''})\ne S(r,f)$. Then by using the similar way as in the proof \cite[Theorem 4]{gs}, we get the equation \eqref{eq.1.3} does not exist exponential polynomial solutions. Here we omit the details.

This completes the proof. $\hfill\Box$

\section{Proof of Theorem \ref{the2}}

Let $f$ be an entire solution of \eqref{eq.1.4}. By \eqref{eq.1.4} and \cite[Theorem 1.22]{f}, we have
\begin{eqnarray*}
    T(r,e^{2d})=m(r,e^{2d})&=&m(r,\frac{af'f''-bf^2}{c})  \\
    & \leqslant& m(r,\frac{af'f''-bf^2}{cf^2})+m(r,f^2)+S(r,f)  \\
     &\leqslant& 2T(r,f)+S(r,f).
\end{eqnarray*}
Then we have $\rho (f)\geqslant k$, where $k=\deg (d) >0$. Taking the logarithm of both sides of the equation \eqref{eq.1.4}, and then take the derivative,
$$\frac{a'f'f''+a(f'')^2+af'f'''-b'f^2-2bff'}{af'f''-bf^2}=\frac{c'+2cd'}{c}.$$
Then
\begin{equation}
\label{4.1}
    ac(f'')^2+acf'f'''-2bcff'-hf'f''=sf^2,
\end{equation}
where $h=ac'-a'c+2acd'$, $s=b'c-bc'-2bcd'$. Noting that $h\not\equiv0$ and $s\not\equiv0$ hold by using the proof of Theorem \ref{the1}. Divide both sides of \eqref{4.1} by $sff'$ to get

$$\frac{ac(f'')^2+acf'f'''-hf'f''}{sff'}-\frac{2bc}{s}=\frac{f}{f'}.$$
Since $a,b,c,d$, and $s$ are polynomials, then from \cite[Theorem 1.22]{f}, we get

$$m(r,\frac{f}{f'})=S(r,f).$$
By differentiating \eqref{4.1}, we obtain
\begin{eqnarray}
\label{4.2}
    (a'c+ac'-h)(f'')^2+(a'c+ac'-h)f'f'''+(4bcd'-4b'c)ff'-  \\ h'f'f''+3acf''f'''+acf'f^{(4)}-2bc(f')^2-2bcff''=s'f^2. \nonumber
\end{eqnarray}
From \eqref{4.1} and \eqref{4.2}, we have
\begin{eqnarray*}
    3acsf''f'''+acsf'f^{(4)}-2bcs(f')^2-2bcsff''+(a'cs+ac's-hs-acs')[(f'')^2+f'f''']  \\+(4bcd's-4csb'+2bcs')ff'+(hs'-h's)f'f''=0.  \nonumber
\end{eqnarray*}
Rewriting as
\begin{eqnarray}
\label{4.3}
    f''[{-2bcsf+(hs'-h's)f'+(a'cs+ac's-hs-acs')f''+3acsf'''}]   \nonumber \\
 =f'[(4csb'-4bcsd'-2bcs')f+2bcsf'-(a'cs+ac's-hs-acs')f'''-acsf^{(4)}].
\end{eqnarray}
Let
\begin{eqnarray}
\label{4.4}
    \zeta&=&\frac{-2bcsf+(hs'-h's)f'+(a'cs+ac's-hs-acs')f''+3acsf'''}{f'}
 \\&=&\frac{(4csb'-4bcsd'-2bcs')f+2bcsf'-(a'cs+ac's-hs-acs')f'''-acsf^{(4)}}{f''}.\nonumber
\end{eqnarray}

By $m(r,\frac{f}{f'})=S(r,f)$ and \eqref{4.4} we have $m(r,\zeta)=S(r,f)$. We claim that the zeros of $f'$ with multiplicity $\ge2$ are at most finitely many, otherwise there are infinitely many of these zeros. Since $ f $ and $f'$ have finitely many common zeros by equation \eqref{eq.1.4}, let $z_0$ is a zero of $f'$ with multiplicity $\ge2$ (not the zero of $f$ and the coefficient of \eqref{4.2}). Substituting $z_0$ into \eqref{4.2}, the left side of \eqref{4.2} is zero and the right side is not zero, this is impossible. Then the zeros of $f'$ with multiplicity $\ge2$ are at most finitely many, which implies that the common zeros of $f'$ and $f''$ are at most finitely many. Now, by \eqref{4.3} and the define of $\zeta$, we get $N(r,\zeta)=O(\log r)$. Then $T(r,\zeta)=m(r,\zeta)+N(r,\zeta)=S(r,f)$.  Equation \eqref{4.4} gives two linear differential equations:
\begin{equation}
\label{4.5}
    3acsf'''+(a'cs+ac's-hs-acs')f''+(hs'-h's-\zeta)f'-2bcsf=0,
\end{equation}
\begin{equation}
\label{4.6}
    -acsf^{(4)}-(a'cs+ac's-hs-acs')f'''-\zeta f''+2bcsf'+(4csb'-4bcsd'-2bcs')f=0.
\end{equation}
By differentiating \eqref{4.5}, we have
\begin{eqnarray}
\label{4.7}
    3acsf^{(4)}+(4a'cs+4asc'+2acs'-hs)f'''  \nonumber  \\
    +(a''cs+2a'c's+asc''-acs''-2h's-\zeta)f'' \nonumber  \\
    +(hs''-h''s-2bcs-\zeta')f'-2(b'cs+bsc'+bcs')f=0.
\end{eqnarray}
From \eqref{4.6} and \eqref{4.7}, we get
\begin{eqnarray}
\label{4.8}
    (a'cs+asc'+5acs'+2hs)f'''+(a''cs+2a'c's+asc''-acs''-2h's-4\zeta)f''  \nonumber \\
    +(hs''-h''s+4bcs-\zeta')f'+(10csb'-8bcs'-2bsc'-12bcsd')f=0.
\end{eqnarray}
From \eqref{4.5} and \eqref{4.8}, we have
\begin{equation}
\label{4.9}
    P_2f''+P_1f'+P_0f=0,
\end{equation}
where
\begin{eqnarray*}
    P_2&=&(a'cs+ac's-hs-acs')(a'cs+asc'+5acs'+2hs)  \\
         &   &   -3acs(a''cs+2a'c's+asc''-acs''-2h's-4\zeta),  \\
    P_1&=&(hs'-h's-\zeta)(a'cs+asc'+5acs'+2hs)-3acs(hs''-h''s+4bcs-\zeta'),  \\
    P_0&=&-2c^2s^2(a'b+15ab')+2bcs(7cs'+2c's)+4bcs^2(9acd'-h).
\end{eqnarray*}
Noting that $P_0\not\equiv0$, otherwise $P_0\equiv0$, divide both sides of $P_0$ by $abc^2s^2$, then
$$-2(\frac{a'}{a}+15\frac{b'}{b})+2(7\frac{s'}{s}+2\frac{c'}{c})=4(\frac{c'}{c}-\frac{a'}{a}+2d')-36d',$$
which implies that $\frac{a^2s^{14}}{b^{30}}=\delta_3e^{-28d}$, where $\delta_3\ne 0$ is a constant, this is a contradiction. Thus $P_0\not\equiv0$.
Differentiating \eqref{4.9}, we have
\begin{equation}
\label{4.10}
    P_2f'''+(P_2'+P_1)f''+(P_1'+P_0)f'+P_0'f=0
\end{equation}
From \eqref{4.5} and \eqref{4.10}, we have
\begin{eqnarray}
\label{4.11}
    [3acs(P_2'+P_1)-P_2(a'cs+ac's-hs-acs')]f''  \nonumber \\
    +[3acs(P_1'+P_0)-P_2(hs'-h's-\zeta)]f'+(3acsP_0'+2P_2bcs)f=0.
\end{eqnarray}
From \eqref{4.9} and \eqref{4.11}, we have
\begin{equation}
\label{4.12}
    uf'=vf,
\end{equation}
where
$$u=[3acsP_2(P_1'+P_0)-P_2^2(hs'-h's-\zeta)]-3acsP_1(P_2'+P_1)+P_1P_2(a'cs+ac's-hs-acs'),$$
$$v=3acsP_0(P_2'+P_1)-P_0P_2(a'cs+ac's-hs-acs')-(3acsP_2P_0'+2P_2^2bcs).$$
Since $a$, $b$, $c$, $d$, $s$ and $h$ are polynomials and $T(r,\zeta)=S(r,f)$, and $\rho(f)\ge k$, then $T(r, P_i)=S(r,f)$, $i=0, 1 ,2$. Therefore, $T(r,u)+T(r,v)=S(r,f)$.

Since $ff'\not\equiv0$, then we can divide the discussion into two cases: $uv\not\equiv 0$ or $u\equiv v\equiv0$.

\indent\textbf{Case 1.}  $uv\not\equiv 0$. From \eqref{4.12}, we get $f'=\frac{vf}{u}$, substituting it into \eqref{eq.1.4}, $$\frac{av}{u}ff''-bf^2=ce^{2d}.$$

We claim that $N(r,\frac{1}{f})=S(r,f)$, otherwise $N(r,\frac{1}{f})\ne S(r,f)$, assume $z_0$ is a zero of $f$ (not the zero of $a$, $b$, $c$, $u$, $v$). Substituting $z_0$ into the above equation, the left side of the above equation is zero and the right side is not zero, this is impossible. Thus $N(r,\frac{1}{f})=S(r,f)$. From \eqref{eq.1.4}, we have
$$\frac{1}{f^2}=\frac{e^{-2d}}{c}(\frac{af'f''}{f^2}-b).$$
From the first fundamental theorem, the above equation and \cite[Theorem 1.22]{f}, we get
\begin{eqnarray*}
    2T(r,f)=T(r,f^2)&=&m(r,\frac{1}{f^2})+N(r,\frac{1}{f^2})+O(1)  \\
    &\leqslant& m(r,\frac{e^{-2d}}{c}(\frac{af'f''}{f^2}-b))+S(r,f)  \\
    &\leqslant& m(r,e^{-2d})+m(r,\frac{f'}{f})+m(r,\frac{f''}{f})+S(r,f)  \\
    &\leqslant& O(r^k)+S(r,f),
\end{eqnarray*}
where $k=\deg (d)$. Thus, $\rho(f)\leqslant k$. Then $\rho(f)= k$. Since $f$ is of finite order, then $N(r,\frac{1}{f})=O(\log r)$. This shows that $f$ has only finitely many zeros, then $f(z)=P_3(z)e^{w(z)}$,  $P_3(z)$ is a polynomial and $w(z)$ is a polynomial of degree $k$.
Substituting $f(z)=P_3(z)e^{w(z)}$ into equation \eqref{eq.1.4}, we get
$${[(aP_3'+aP_3w')(P_3''+2P_3'w'+P_3w''+P_3(w')^2)-bP_3^2]}e^{2w}=ce^{2d}.$$
There's only one constant difference between $w(z)$ and $d(z)$. Therefore, we could just write $f(z)=P(z)e^{d(z)}$, where $P(z)$ is a polynomial such that
$$(aP'+aPd')[P''+2P'd'+Pd''+P(d')^2]-bP^2\equiv c,$$
where $c$ is a non-constant polynomial. This is the conclusion  ${\textnormal{(i$_1$)}}$ of Theorem \ref{the2}.

Next, if $c$ is a non-zero constant, then we have
\begin{equation}
\label{4.13}
    aP'P''+2a(P')^2d'+aPP'd''+3aPP'(d')^2+aPP''d'+P^2(ad'd''+a(d')^3-b)\equiv c.
\end{equation}

\indent\textbf{Subcase 1.1.}    $\deg (d)=1$, set $d(z)=\lambda z +c_3$, where $\lambda\ne 0$ and $c_3$ are constants. Then \eqref{4.13} turns into
\begin{equation}
\label{4.14}
    aP'P''+2\lambda a(P')^2+3a\lambda ^2PP'+\lambda aPP''+P^2(\lambda^3a-b)\equiv c.
\end{equation}

\indent\textbf{Subcase 1.1.1.} Both $a$ and $b$ are constants. If $P$ is a non-constant polynomial, then the highest degree term of equation \eqref{4.14} is $P^2(\lambda^3a-b)$, and the right side of \eqref{4.14} is a constant, which is a contradiction. Then $\lambda^3a-b=0$,  $a=\frac{b}{\lambda^3}$, put it into  equation \eqref{4.14} we have
$$\frac{b}{\lambda^3}P'P''+\frac{2b}{\lambda^2}(P')^2+\frac{3b}{\lambda}PP'+\frac{b}{\lambda^2}PP''\equiv c.$$

As can be seen from the above equation, the degree of $\frac{3b}{\lambda}PP'$ is higher than other terms, which is a contradiction, then $\frac{3b}{\lambda}PP'$ is a constant, which means that $b$, $P$ are constants. This contradicts with the assumption that $P$ is a polynomial. Thus $P$ is a non-zero constant. Then from
\eqref{4.14}, we have $P^2(\lambda^3a-b)=c$, where $\lambda^3a-b$ is a non-zero constant. Then the solution $f$ is of form $f(z)=c_1e^{\lambda z}$, where $c_1=Pe^{c_3}$, and $c_1^2(\lambda^3a-b)=ce^{2c_3}$ by $P^2(\lambda^3a-b)=c$. This is the conclusion ${\textnormal{(ii$_1$)}}$ of Theorem \ref{the2}.

\indent\textbf{Subcase 1.1.2.} There exists at least one of $a$ and $b$ is not constant.  If $P$ is a non-constant polynomial, we claim that $\deg (a) =\deg (b)$, otherwise, without loss of generality, assume $\deg (a) >\deg (b)$, then the highest degree term of equation \eqref{4.14} is $P^2(\lambda^3a-b)$, and the right side of \eqref{4.14} is a constant, which is a contradiction. Thus $\deg (a) =\deg (b)$. If $P$ is a non-zero constant, we claim that $\deg (a) =\deg (b)$, otherwise, without loss of generality, assume $\deg (a) >\deg (b)$, then from \eqref{4.14}, we have $P^2(\lambda^3a-b)=c$, this is a contradiction. Thus $\deg (a) =\deg (b)$. Then the solution $f$ is of form $f(z)=Pe^{d(z)}$, where $P$ is a non-zero constant or polynomial. This is the conclusion ${\textnormal{(i$_{21}$)}}$ of Theorem \ref{the2}.

\indent\textbf{Subcase 1.2.}  $\deg (d)\ge 2$.

\indent\textbf{Subcase 1.2.1.}   $ad'd''+a(d')^3-b\equiv0$.

If $P$ is a non-zero constant, then from \eqref{4.13}, we get $c=0$, which is a contradiction. Thus $P$ is a non-constant polynomial. Then $3aPP'(d')^2$ has higher degree than the other terms, and the right side of \eqref{4.13} is a constant, which is a contradiction. Thus the case for  $ad'd''+a(d')^3-b\equiv0$ does not exist.

\indent\textbf{Subcase 1.2.2.}  $ad'd''+a(d')^3-b= M_1$, where $M_1$ is a non-zero constant.

We claim $P$ is a non-zero constant, otherwise, from \eqref{4.13}, $3aPP'(d')^2$ has higher degree than the other terms, and the right side of \eqref{4.13} is a constant, which is a contradiction. Thus $P$ is constant. From \eqref{4.12}, we have $P^2M_1=c$. Then the solution $f$ is of form $f(z)=Pe^{d(z)}$. This is the conclusion ${\textnormal{(i$_{22}$)}}-1)$ of Theorem \ref{the2}.

\indent\textbf{Subcase 1.2.3.} $ad'd''+a(d')^3-b=M_2$, where $M_2$ is a polynomial with $\deg (M_2)\ge 1$. Then \eqref{4.13} becomes
\begin{equation}
\label{4.15}
    aP'P''+2a(P')^2d'+aPP'd''+3aPP'(d')^2+aPP''d'+P^2M_2= c.
\end{equation}

We claim that $P$ is non-constant polynomial, otherwise from \eqref{4.15} we get $P^2M_2=c$, this is impossible. Observing \eqref{4.15}, we claim that $\deg (P^2M_2)=\deg (3aPP'(d')^2)$, otherwise, without loss of generality, suppose $\deg (P^2M_2)>\deg (3aPP'(d')^2)$, from \eqref{4.15} we get a contradiction. Thus $\deg (P^2M_2)=\deg (3aPP'(d')^2)$, i.e. $\deg (M_2)=\deg (a) +2\deg (d)-3$. This is the conclusion ${\textnormal{(i$_{22}$)}}-2)$ of Theorem \ref{the2}.

\indent\textbf{Case 2.} $u\equiv v\equiv0.$

We claim that $P_2\not\equiv0$, otherwise from $v\equiv0$, we get $3acsP_0P_1\equiv0$, then we obtain $P_1\equiv0$, and then \eqref{4.9} becomes $P_0f=0$, which is a contradiction. Thus $P_2\not\equiv0$. Rewriting $P_2,P_1$ as
$$P_2=A_0\zeta+A_1,   $$
where $A_0=12acs\not\equiv0$, $A_1=(a'cs+ac's-hs-acs')(a'cs+asc'+5acs'+2hs)-3acs(a''cs+2a'c's+asc''-acs''-2h's)$ is a polynomial,
$$P_1=B_0\zeta+B_1\zeta'+B_2,$$
where $B_0=-ac's-asc'-5acs'-2hs$, $B_1=3acs$, $B_2$ is a polynomial. By $v\equiv0$, we have
\begin{equation}
\label{4.16}
    -2bcsP_2^2+[-3acsP_0'-(a'cs+ac's-hs-acs')P_0]P_2+3acsP_0(P_2'+P_1)\equiv0.
\end{equation}
One can simplify \eqref{4.16} into new form, that is
\begin{equation}
\label{4.17}
    X_1P_2^2+X_2P_2+X_3(P_2'+P_1)\equiv0,
\end{equation}
where $X_1=-2bcs$, $X_2=-3acsP_0'-(a'cs+ac's-hs-acs')P_0$, $X_3=3acsP_0$. Substituting $P_2$ and $P_1$ into \eqref{4.17}, we have
\begin{eqnarray}
\label{4.18}
    X_1A_0^2\zeta^2+(2A_0A_1X_1+X_2A_0+X_3A_0'+X_3B_0)\zeta  \nonumber \\
    +X_3(A_0+B_1)\zeta'+(X_1A_1^2+X_2A_1+X_3A_1'+X_3B_2)\equiv0.
\end{eqnarray}

 We claim $\zeta$ is a rational function. Otherwise,
 suppose $\zeta$  is a transcendental meromorphic function. Since $X_1A_0^2\not\equiv0$, using Clunie Lemma \cite[Lemma 3.3]{c} to \eqref{4.18}, then $m(r,\zeta)=S(r,\zeta)$. Noting that  $N(r,\zeta)=O(\log r)$, then $T(r,\zeta)=O(\log r)+S(r,\zeta)$, which is impossible.

If $s$ and $a$ are both non-zero constants, then by comparing the degree of  both sides of the equation $s=b'c-bc'-2bcd'$,  then $2bcd'$ must be constant. Hence, $a$, $b$, $c$, $d'$ and $s$ are non-zero constants, then we can get $A_0$, $A_1$, $B_i$, $X_i$, $(i=1, 2, 3 )$ are constants.

Now we claim $\zeta$ has no poles. Otherwise, suppose $z_0$ is a pole of $\zeta$ with multiplicity $n$. Substituting $z_0$ into \eqref{4.18}, we see that $n=1$, i.e. $z_0$ is a simple pole. Therefore $\zeta=\frac{l}{z-z_0}+O(1)$, as $z\rightarrow z_0$, where $l$ is a nonzero constant.  Substituting $\zeta=\frac{l}{z-z_0}+O(1)$ into $u=3acsA_0B_1\zeta\zeta''-3acsB_1(A_0+B_1)(\zeta')^2+P(z,\zeta)\equiv0$, where $P(z,\zeta)$ is a differential polynomial of $\zeta$ and the weight of $P(z,\zeta)$ is less than or equal to $3$. Then we must have $3A_0-B_1=33acs\equiv0$
, which is a contradiction. Thus $\zeta$ has no poles, which means $\zeta$  is a polynomial. By comparing the degree of  both sides of  \eqref{4.18}, we see $\zeta$ is a constant. Then \eqref{4.9} becomes
\begin{equation}
\label{4.19}
    f''+\frac{P_1}{P_2}f'+\frac{P_0}{P_2}f=0,
\end{equation}
where $P_0,P_1,P_2$ are non-zero constants. Then the characteristic equation of  \eqref{4.19} is
$$t^2+\frac{P_1}{P_2}t+\frac{P_0}{P_2}=0,$$
where $t_1$ and $t_2$ are two characteristic roots. Since $P_0\not \equiv0$, then $t_1$ and $t_2$ are non-zero constants.

\indent\textbf{Subcase 2.1.} $t_1\ne t_2$. Then
\begin{equation}
\label{4.20}
    f(z)=c_1e^{t_1z}+c_2e^{t_2z},
\end{equation}
where $c_1$ and $c_2$ are constants. Next we will discuss the case of whether $c_1$, $c_2$ is zero.

\indent\textbf{Subcase 2.1.1.} $c_1c_2=0$. Without loss of generality, assume that $c_2=0$.  Since $d'$ is a non-zero constant, then assume $d=\lambda z+c_3$, where $\lambda (\ne 0)$ and $c_3$ are constants. Substituting $f(z)=c_1e^{t_1z}$ and  $d=\lambda z+c_3$ into equation \eqref{eq.1.4} to get
$$c_1^2(at_1^3-b)e^{2t_1z}=ce^{2c_3}e^{2\lambda z},$$
which implies that $\lambda=t_1$, and $c_1^2(a\lambda^3-b)=ce^{2c_3}$. Then the solution $f$ is of form $f(z)=c_1e^{\lambda z}$. This is the conclusion  ${\textnormal{(ii$_1$)}}$ of Theorem \ref{the2}.

\indent \textbf{Subcase 2.1.2.}  $c_1c_2\ne0$. Substituting \eqref{4.20} and $d=\lambda z+c_3$ into equation \eqref{eq.1.4}, we have
\begin{equation}
\label{4.21}
    c_1^2(a t_1^3-b)e^{2t_1z}+c_2^2(at_2^3-b)e^{2t_2z}+(a c_1c_2t_1t_2^2+a c_1c_2t_1^2t_2-2b c_1c_2)e^{(t_1+t_2)z}=ce^{2c_3}e^{2\lambda z}.
\end{equation}
Applying \cite[Theorem 1.51]{f} to equation \eqref{4.21}, then $\lambda= t_1$ or $\lambda=t_2$ or $\lambda=\frac{t_1+t_2}{2}$.

If $\lambda=t_1$, then \eqref{4.21} can be written as
$$[c_1^2(at_1^3-b)-ce^{2c_3}]e^{2t_1z}+c_2^2(a t_2^3-b)e^{2t_2z}+(a c_1c_2t_1t_2^2+a c_1c_2t_1^2t_2-2b c_1c_2)e^{(t_1+t_2)z}=0.$$
By applying \cite[Theorem 1.51]{f} to the above equation, we have
$$c_1^2(at_1^3-b)-ce^{2c_3}=c_2^2(at_2^3-b)=a c_1c_2t_1t_2^2+a c_1c_2t_1^2t_2-2b c_1c_2=0,$$
which gives $at_2^3=b$, and $c_1c_2t_2(2t_2+t_1)(t_1-t_2)=0$. Since $t_1\ne t_2$, then we have $\lambda=t_1=-2t_2$. Combining this with the equation  \eqref{4.20}, we have $f(z)=c_1e^{\lambda z}+c_2e^{-\frac{\lambda}{2}z}$ and $\frac{9}{8}c_1^2a
\lambda^3=ce^{2c_3}$. This is the conclusion ${\textnormal{(ii$_2$)}}$ of Theorem \ref{the2}.

If $\lambda=t_2$, using the same method as $\lambda=t_1$, we get $\lambda=t_2=-2t_1$. Then from \eqref{4.20},  we have $f(z)=c_1e^{-\frac{\lambda}{2} z}+c_2e^{\lambda z}$ and $\frac{9}{8}c_2^2a\lambda^3=ce^{2c_3}$. This is the conclusion ${\textnormal{(ii$_2$)}}$ of Theorem \ref{the2}.

If $\lambda=\frac{t_1+t_2}{2}$, then \eqref{4.21} becomes
$$c_1^2(a t_1^3-b)e^{2t_1z}+c_2^2(a t_2^3-b)e^{2t_2z}+(ac_1c_2t_1t_2^2+ac_1c_2t_1^2t_2-2b c_1c_2-ce^{2c_3})e^{(t_1+t_2)z}=0.$$
By applying \cite[Theorem 1.51]{f} to the above equation, we have
$$c_1^2(a t_1^3-b)=c_2^2(a t_2^3-b)=a c_1c_2t_1t_2^2+a c_1c_2t_1^2t_2-2b c_1c_2-ce^{2c_3}=0.$$
Then we get $t_1^3=t_2^3=\frac{b}{a}$ and $c_1c_2t_1(t_2-t_1)(t_2+2t_1)=ce^{2c_3}$. From \eqref{4.20}, we have $f(z)=c_1e^{t_1z}+c_2e^{t_2z}$. This is the conclusion ${\textnormal{(ii$_3$)}}$ of Theorem \ref{the2}.

\indent\textbf{Subcase 2.2.} $t_1= t_2$. Then
\begin{equation}
\label{4.22}
    f(z)=(c_4+c_5 z)e^{t_1 z},
\end{equation}
where $c_4, c_5$ are constants. Substituting \eqref{4.22} into equation \eqref{eq.1.4}, we obtain
$$[a(c_5+c_4t_1+c_5t_1z)(2c_5t_1+c_4t_1^2+c_5t_1^2z)-b(c_4^2+2c_4c_5z+c_5^2z^2)]e^{2t_1z}=ce^{2c_3} e^{2\lambda z},$$
which implies that $\lambda=t_1$ and
\begin{equation}
\label{4.23}
    a(c_5+c_4t_1+c_5t_1z)(2c_5t_1+c_4t_1^2+c_5t_1^2z)-b(c_4^2+2c_4c_5z+c_5^2z^2)=ce^{2c_3}.
\end{equation}
By considering the coefficients at $z^2$ and $z$, respectively, we have
$$c_5^2(a t_1^3-b)=0,$$
and$$3a c_5^2t_1^2+2a c_4c_5t_1^3-2b c_4c_5=0.$$
If $c_5\ne 0$, then $at_1^3=b$, and $3ac_5^2t_1^3=0$, which implies that $c_5=0$, this is a contradiction. Thus $c_5=0$. Substituting $c_5=0$ into \eqref{4.23}, we obtain $c_4^2(a\lambda^3-b)=ce^{2c_3}$. Then the solution $f$ is of form $f(z)=c_4e^{\lambda z}$. This is the conclusion ${\textnormal{(ii$_1$)}}$ of Theorem \ref{the2}.

If there exists at least one of $s$ and $a$ is a non-constant polynomial, then by \cite[Remark 1 of Theorem 4.1]{g}, the entire solution of the equation \eqref{4.9} is of finite order. Thus $\lambda(f)\leqslant \rho(f)<\infty$. This is the conclusion ${\textnormal{(i$_3$)}}$ of Theorem \ref{the2}.

The proof of Theorem \ref{the2} is now completed.  $\hfill\Box$


\end{document}